\documentclass[11pt]{article}
\usepackage{pstricks,psfrag,theorem,cite,enumerate}
\usepackage{times,graphics,epsfig}
\usepackage{amssymb}
\usepackage{amsmath}
\usepackage{amsfonts}
\usepackage{amsfonts,amsmath,amstext}
\usepackage{psfrag,latexsym,verbatim,graphicx}
\usepackage{latexsym}
\usepackage{colordvi}
\usepackage{epsfig}
\usepackage{graphicx}

\parskip\smallskipamount

\newtheorem{theorem}{Theorem}

\newtheorem{lem}{Lemma}

\newtheorem{cor}{Corollary}

\newtheorem{prop}{Proposition}

\newtheorem{rem}{Remark}

\newcommand{\R}{\mathbb{R}}

\newcommand{\Z}{\mathbb{Z}}
\newcommand{\W}{\mathbb{W}}
\newcommand{\N}{\mathbb{N}}
\newcommand{\Q}{\mathbb{Q}}
\newcommand{\HH}{\mathbb{H}}
\newcommand{\D}{\mathbb{D}}

\renewcommand{\phi}{\varphi}
\renewcommand{\epsilon}{\varepsilon}

\parskip \smallskipamount

\title{Poisson Hail
on a Hot Ground}

\author{Francois Baccelli\footnote{
\emph{INRIA and ENS, Paris, France}. The research of this author was partially supported by EURONF.}
\
 and Sergey Foss\footnote{
\emph{Heriot-Watt University, Edinburgh, UK and Institute of Mathematics,
Novosibirsk, Russia}. The research of this author was partially supported by
INRIA and EURONF.\hfil\hspace{11cm}\hfil
The authors thank the SCS Programme of Isaac Newton Institute for Mathematical Studies
where this work was completed.}}

\date{13 March 2011}

\begin{document}

\maketitle

\begin{abstract}

We consider a queue where the server is the Euclidean space, and
the customers are
random closed sets (RACS) of the Euclidean space. These RACS arrive according
to a Poisson rain and each of them has a random service time
(in the case of hail falling on the Euclidean plane, this is the height of the hailstone,
whereas the RACS is its footprint). The Euclidean space serves
customers at speed 1.  The service discipline is a hard exclusion rule:
no two intersecting RACS can be served simultaneously and service is in
the First In First Out order: only the hailstones in contact with the ground melt at speed 1,
whereas the other ones are queued; a tagged RACS waits until all RACS arrived before
it and intersecting it have fully melted before starting its own melting. We give the
evolution equations for this queue.
We prove that it is stable for a sufficiently small arrival intensity, provided the typical
diameter of the RACS and the typical service time have finite exponential moments. We also
discuss the percolation properties of the stationary regime of the RACS in the queue.

\end{abstract}

\begin{quotation}\small

\end{quotation}

{\bf Keywords:}

Poisson point process, Poisson rain, random closed sets, Euclidean space,
service, stability, backward scheme, monotonicity, branching process, percolation,
hard core exclusion processes, queueing theory, stochastic geometry.

\section{Introduction}

Consider a Poisson rain on the $d$ dimensional Euclidean space $\R^d$ with intensity $\lambda$; by Poisson rain,
we mean a Poisson point process of intensity $\lambda$ in $\R^{d+1}$
which gives the (random) number of arrivals in all time-space Borel sets.
Each Poisson arrival, say at
location $x$ and time $t$, brings a customer with two main characteristics:

\begin{itemize}
\item A grain $C$, which is a RACS of
  $\R^d$ \cite{Stoyan:1995:SGA} {\em centered} at the origin. 
If the RACS is a ball with random radius, its center is that of the ball.
For more general cases, the center of a RACS could be defined as e.g. its gravity center.
\item A random service time $\sigma$.
\end{itemize}
In the most general setting, these two characteristics will be assumed to be marks of the point process.
In this paper, we will concentrate on the simplest case, which is that of an independent marking and 
independent and identically distributed (i.i.d.) marks: the mark
$(C,\sigma )$ of point $(x,t)$ has 
 some given distribution and is independent of everything else.

The customer arriving at time $t$ and location $x$ with mark $(C,\sigma)$ creates a hailstone, with footprint $x+C$ in $\R^d$ 
and with height $\sigma$.

These hailstones do not move: they are to be melted/served by the Euclidean plane at the
location where they arrive in the FCFS order, respecting some
hard exclusion rules:
if the footprints of two hailstones have a non empty intersection,
then the one arriving second has to wait for the end of the melting/service of the
first to start its melting/service.
Once the service of a customer is started, it proceeds uninterrupted at speed 1. Once a customer
is served/hailstone fully melted, it leaves the Euclidean space.

Notice that the customers being served at any given time form a
hard exclusion process as no two customers having intersecting footprints are ever served
at the same time. For instance, if the grains are balls, the footprint
balls concurrently served form a hard ball exclusion process.
Here are a few basic questions on this model:

\begin{itemize}
\item Does there exist any positive
$\lambda$ for which this model is (globally) stable?
By stability, we mean that, for all $k$ and for all
bounded Borel set $B_1,\ldots,B_k$, the vector
$N_1(t),\ldots,N_k(t)$, where $N_j(t)$ denotes the number of RACS
which are queued or in service
at time $t$ and intersect the Borel set $B_j$, converges in distribution to a finite
random vector when $t$ tends to infinity.
\item If so, does the stationary regime
percolate? By this, we mean that the union of the RACS which are queued or in service in
a snapshot of the stationary regime has an infinite connected component.
\end{itemize}

The paper is structured as follows.
In section \ref{GM}, we study pure growth models (the ground is cold and hailstones do not melt)
and show that the heap formed by the customers grows with
(at most) linear rate with time and that the growth rate
tends to zero if the input rate tends to zero.
We consider models with service (hot ground)
in section \ref{secservar}.
Discrete versions of the problems are studied in section \ref{secbhhg}.

\section{Main Result}

Our main result bears on the construction of the stationary regime of this system.

As we shall see below (see in particular Equations (\ref{eq:hail1}) and (\ref{eq:simpledyn})),
the Poisson Hail model falls in the category of infinite dimensional max plus linear systems.
This model has nice monotonicity properties (see sections \ref{GM} and \ref{secservar}).
However it does not satisfy the separability property of \cite{sep}, which
prevents the use of general sub-additive ergodic theory tools to assess stability,
and makes the problem interesting.

Denote by $\xi$ the (random) diameter of the typical RACS (i.e. the maximal distance between
its points) and by $\sigma$ the service time of that RAC. Assume that the system starts at
time $t=0$ from the empty state and denote by $W^x_t$ the time to empty the system of
all RACS that contain point $x$ and that arrive by time $t$.

\begin{theorem}\label{thmain}
Assume that the Poisson hail starts at time $t=0$ and that the
system is empty at that time. Assume further that the
distributions of the random variables $\xi^d$ and $\sigma$ are
{\it light-tailed}, i.e. there is a positive constant $c$ such
that ${\mathbf E} e^{c\xi^d}$ and ${\mathbf E} e^{c\sigma}$ are
finite. Then there exists a positive constant $\lambda_0$ (which
depends on $d$ and on the joint distribution of $\xi$ and
$\sigma$) such that, for any $\lambda < \lambda _0$, the model is
globally stable. This means that, for any finite set $A$ in
$\R^d$, as $t\to\infty$, the distribution of the random field
$(W^x_t, \ x\in A)$ converges weakly to the stationary one.
\end{theorem}

\section{Growth Models}
\label{GM}
Let $\Phi$ be a marked Poisson point process in $\R^{d+1}$:
for all Borel sets $B$ of $\R^d$ and $a\le b$, a r.v.
$\Phi(B,[a,b])$
denotes the number of RACS with center
located in $B$
that arrive in the time interval $[a,b]$. The marks of this point process
are i.i.d. pairs $(C_n,\sigma_n)$, where $C_n$ is a RACS of $\R^d$ and $\sigma_n$
is a height (in $\R+$, the positive real line).

The growth model is best defined by the following equations
satisfied by $H_t^x$, the height at location $x\in \R^d$
of the heap made of all RACS arrived before time
$t$ (i.e. in the $(0,t)$ interval): for all $t>u\ge 0$,
\begin{equation}\label{eq:hail1}
H_t^x = H_u^x + \int_{[u,t)}  \left(\sigma_v^x +\sup_{y\in C_v^x } H^y_v  -H_v^x\right) N^x(dv),
\end{equation}
where $N^x$ denotes the Poisson point process on $\R^+$ of RACS arrivals
intersecting location $x$:
$$ N^x ([a,b]) = \int_{\R^d \times [a,b]} 1_{C_v \cap \{x\}\ne \emptyset} \Phi(dv),$$
and $\sigma_u^x$ (resp. $C^x_u$) the canonical height (resp. RAC) mark process of $N^x$.
That is, if the point process $N^x$ has points $T^x_i$, and if one denotes by
$(\sigma_i^x,C^x_i)$ the mark of point $T^x_i$, then $\sigma_u^x$ (resp. $C^x_u$) is
equal to $\sigma_i^x$ (resp. $C^x_i)$) on $[T^x_i,T^x_{i+1})$.

These equations lead to some measurability questions. Below, we will
assume that the RACS are such that the last supremum actually bears on a subset of $\Q^d$,
where $\Q$ denotes the set of rational numbers,
so that these questions do not occur.

Of course, in order to specify the dynamics, one also needs some initial condition, namely
some initial field $H^x_0$, with $H^x_0\in \R$ for all $x\in \R^d$.

If one denotes by $\tau^x(t)$ the last epoch of $N^x$ in $(-\infty,t)$,
then this equation can be rewritten as
the following recursion:
\begin{eqnarray*}
H_t^x & \hspace{-.2cm}= \hspace{-.2cm}& H_0^x + \int_{[0,\tau^x(t))}
 \left(\sigma_v^x +\sup_{y\in C_v^x } H^y_v  -H_v^x\right) N^x(dv)
+ \sigma_{\tau^x(t)}^x +
\sup_{y\in C^x_{\tau^x(t)}} H_{\tau^x(t)}^y - H_{\tau^x(t)}^x~,
\end{eqnarray*}
that is
\begin{equation}\label{eq:simple2}
H_t^x
=
\left(\sigma_{\tau^x(t)}^x +
\sup_{y\in C^x_{\tau^x(t)}} H_{\tau^x(t)}^y\right)
1_{\tau^x(t)\ge 0}
+H_{0}^x1_{\tau^x(t)<0}.
\end{equation}
These are the forward equations. We will also use the backward equations,
which give the heights at time $0$ for an arrival point process
which is the restriction of the Poisson hail to the interval $[-t,0]$ for $t>0$.
Let $\HH^x_t$ denote the height at locations $x$ and time $0$ for this point process.
Assuming that the initial condition is 0, we have
\begin{equation}\label{eq:simple2b}
\HH_t^x = \left(\sigma_{\tau^x_-(t)}^x +
\sup_{y\in C^x_{\tau^x_-(t)}} \HH_{t+\tau^x_-(t)}^y\circ \theta_{\tau^x_-(t)}\right)
1_{\tau^x_-(t)\ge -t},
\end{equation}
with
$\tau^x_-(t)$ the last arrival of the point process $N^x$ in the interval $[-t,0]$, $t>0$,
and with $\{\theta_u\}$ the time shift on the point processes \cite{BB}.
\begin{rem}
Here are a few important remarks on these Poisson hail equations:
\begin{itemize}
\item The last pathwise equations hold for all point processes
and all RACS/heights (although one has to specify how to handle ties when 
RACS with non-empty intersection arrive at the same time - we postpone the
discussion on this matter to section \ref{secbhhg}).
\item These equations can be extended to the case where customers have a more general
structure than the product of a RACS of $\R^d$ and an interval of the form $[0,\sigma]$. We will call as 
{\em profile} a function $s(y,x):\R^d\times \R^d \to \R\cup \{-\infty\}$, where $s(y,x)$
gives the {\em height} at $x$ relative to a point $y$; we will say 
that point $x$
is constrained by point $y$ in the profile if $s(y,x)\ne -\infty$. The equations
for the case where random profiles (rather than product form RACS) arrive are
\begin{equation}\label{eq:profile}
H_t^x
=
\left(
\sup_{y\in \R^d}
\left(
H_{\tau^x(t)}^y
+
s_{\tau^x(t)}(y,x)
\right)\right)
1_{\tau^x(t)\ge 0}
+H_{0}^x1_{\tau^x(t)<0},
\end{equation}
where $\tau^x(t)$ is the last date of arrival of $N^x$ before time $t$, with $N^x$ the point
process of arrivals of profiles having a point which constrains $x$.
We assume here that this point process has a finite intensity.
The case of product form RACS considered above is a special case with
$$ s_{\tau^x(t)}(y,x)= \begin{cases}
\sigma_{\tau^x(t)} & \mbox{if}\ y\in C^x_{\tau^x(t)}\\
-\infty & \mbox{otherwise}\end{cases},$$
with $N^x$ the point process of arrivals with RACS intersecting $x$.
\end{itemize}
\end{rem}
Here are now some monotonicity properties of these equations:
\begin{enumerate}
\item The representation \eqref{eq:simple2} shows that if we have two marked point processes
$\{N^x\}_x$ and $\{\widetilde N^x\}_x$ such that for all $x$, $N^x\subset \widetilde N^x$
(in the sense that each point of $N^x$ is also a point of $\widetilde N^x$), and if the
marks of the common points are unchanged, then
$H_{t}^x\le \widetilde H_t^x$ for all $t$ and $x$ whenever $H_{0}^x\le \widetilde H_0^x$ for all $x$.
\item Similarly, if we have two marked point processes
$\{N^x\}_x$ and $\{\widetilde N^x\}_x$ such that for all $x$,
$N^x\le \widetilde N^x$ (in the sense that for all $n$, the $n$-th point of $N^x$ is later than the
$n$-th point of $\widetilde N^x$), and the marks are unchanged, then
$H_{t}^x\le \widetilde H_t^x$ for all $t$ and $x$ whenever $H_{0}^x\le \widetilde H_0^x$ for all $x$.
\item Finally, if the marks of a point process are changed in such a way that
$C\subset \widetilde C$ and $\sigma\le \widetilde \sigma$,
then
$H_{t}^x\le \widetilde H_t^x$ for all $t$ and $x$ whenever $H_{0}^x\le \widetilde H_0^x$ for all $x$.
\end{enumerate}
These monotonicity properties hold for the backward construction as well.

They are also easily extended to profiles.
For instance, for the last monotonicity property, if profiles are changed in such a way that
$$ s(y,x)\le \widetilde  s(y,x),\quad \forall x,y,$$
then $H_{t}^x\le \widetilde H_t^x$ for all $t$ and $x$ whenever $H_{0}^x\le \widetilde H_0^x$ for all $x$.

Below, we use these monotonicity properties to get upper-bounds on the $H^x_t$ and $\HH^x_t$ variables.

\subsection{Discretization of Space}\label{ssds}
Consider the lattice $\Z^d$, where $\Z$ denotes the set of integers.
To each point in $x=(x_1,\ldots ,x_d)\in \R^d$, we associate
the point $z(x)=(z_1(x),\ldots ,z_d(x))\in \Z^d$ with coordinates $z_i(x) = \lfloor x_i\rfloor$ where
$\lfloor \cdot \rfloor$ denotes the integer-part. Then, with
the RACS $A$ centered at point $x\in \R^d$ and having diameter $\xi$, we associate an auxiliary
RACS $\breve{A}$ centered at point $z(x)$ and being the $d$-dimensional cube of
side $2\lfloor \xi\rfloor +2$. Since $A \subseteq \breve{A}$,
when replacing the RACS $ A$ by
the RACS $\breve{A}$ at each arrival, and keeping all other features
unchanged, we get from the monotonicity property 3 that for all $t\in \R$ and $x\in \R^d$,
$$
H_t^x \le \breve H_t^{z(x)},
$$
with $\breve H_t^z$ the solution of the discrete state space recursion
\begin{equation}
\label{eq:discrete}
\breve H_t^z = \left(\sigma_{\breve \tau^z(t)}^z +
\max _{y\in \Z^d\cap\breve C_{\breve \tau^z(t)}^z}
\breve H_{\breve \tau^z(t)}^y\right)
1_{\breve \tau^z(t)\ge 0}
+\breve{H}_{0}^z1_{\breve \tau^z(t)<0}
,\quad z\in \Z^d~,
\end{equation}
with $\breve \tau^z(t)$ the last epoch of the point process
$$\breve N^z ([a,b]) = \int_{\R^d \times [a,b]} 1_{\breve C_v \cap \{z\}\ne \emptyset} \Phi(dv)$$
in $(-\infty,t)$.
The last model will be referred to as Model 2.
We will denote by $R$ the typical half-side of the cubic RACS in this model.
These sides are i.i.d. (w.r.t. RACS), and if $\xi^d$ has a light-tailed distribution, then
$R^d$ has too.

\subsection{Discretization of Time}\label{ssdistim}

The discretization of time is in three steps.\\

\noindent{\em Step 1.}
Model 3 is defined as follows: all RACS centered on $z$ that arrive to Model 2
within time interval $[n-1,n)$, arrive to Model 3 at time instant $n-1$.
The ties are then solved according to the initial continuous time ordering. In view of the
monotonicity property 2, Model 3 is an upper bound to Model 2.

Notice that for each $n$, the
arrival process at time $n$ forms a discrete Poisson field of parameter $\lambda$, i.e.
the random number of RACS $M^z_n$ arriving at point $z\in \Z^d$ at time $n$ has a Poisson distribution
with parameter $\lambda$, and these random variables are i.i.d. in $z$ and $n$.

Let $(R_{n,i}^z,\sigma_{n,i}^z)$, $i=1,2,\ldots,M_n^z$, be the i.i.d. radii and heights of the cubic
RACS arriving at point $z$ and time $n$. Let further $M=M_0^0$, $R_i=R_{0,i}^0$, and
$\sigma_i=\sigma_{0,i}^0$.\\

\noindent{\em Step 2.}
Let $R^{z, max}_n$ be the maximal half-side of all RACS that arrive at point $z$ and
time $n$ in Model 2, and $R^{max}=R^{z, max}_n$.
The random variables $R^{z,max}_n$ are i.i.d. in $z$ and in $n$.
We adopt the convention that $R^{z,max}_n=0$ if there is
no arrival at this point and this time.
If the random variable $\xi^d$ is light-tailed, the distribution of $R^d$ is also light-tailed,
and so is that of $\left(R^{max}\right)^d$. Indeed,
$$
\left(R^{max}\right)^d = \left(\max_{i=1}^M R_{i}\right)^d \le \sum_1^M R_{i}^d,
$$
so, for $c >0$,
$$
{\mathbf E} e^{c \left(R^{max}\right)^d} \le
{\mathbf E} e^{c \sum_1^M R_i^d}
=\exp \left(\lambda {\mathbf E} e^{c R^d}\right) < \infty
$$
given ${\mathbf E} e^{c R^d}$ is finite.
Let
$$
\sigma_n^{z,sum}= \sum_{i=1}^{M_n^z} \sigma_{n,i}^z \quad \mbox{and}
\quad
\sigma^{sum} = \sigma_0^{0,sum}.
$$
Then, by
Similar arguments,
$
\sigma^{sum}
$
has a light-tailed distribution if $\sigma_i$ do.
By monotonicity property 3 (applied to the profile case), when replacing the heap of RACS
arriving at $(z,n)$ in Model 3 by the cube of half-side $R^{z,max}_n$ and
of height $\sigma^{z,sum}_n$, for all $z$ and $n$, one again gets an upper bound system
which will be referred to as Model 4.\\

\noindent{\em Step 3.}
The main new feature of the last discrete time Models (3 and 4) is that the RACS that arrive at
some discrete time on different sites may overlap. Below,
we consider the clump made by overlapping RACS as a profile and use monotonicity
property 3 to get a new upper bound model, which will be referred to as Boolean Model 5.

Consider the following discrete Boolean model, associated with time $n$. We say that there is
a "ball" at $z$ at time $n$ if $M_n^z \ge 1$ and that there is no ball at $z$ at this time otherwise.
By ball, we mean a $L_{\infty}$ ball with center $z$ and radius $R^{z,max}_n$.
By decreasing $\lambda$, we can make the probability
$p = {\mathbf P} (M_n^z\ge 1)$ as small as we wish.

Let $\widehat C_n^z$ be the {\em clump} containing point $z$ at time $n$,
which is formally defined as follows:
if there is a ball at $(z,n)$, or another ball of time $n$ covering $z$, this clump is
the largest union of connected balls (these balls are considered as subsets of $\Z^d$ here)
which contains this ball at time $n$; otherwise, the clump is {\em empty}.
For all sets $A$ of the Euclidean space, let $L(A)$ denote the
number of points of the lattice $\Z^d$ contained in $A$.
It is known from percolation theory that, for $p$ sufficiently small,
this clump is a.s. finite \cite{hall} and, moreover, $L(\widehat C_n^z)$
has a light-tailed distribution (since $\left(R^{max}\right)^d$ is light-tailed) \cite{BRS}.
Recall that the latter means that ${\mathbf E} \exp (c L(\widehat C_0^z)) < \infty$, for some $c >0$.

Below, we will denote by $\lambda_c$ the critical value of $\lambda$ below which
this clump is a.s. finite and light-tailed.

For each clump $\widehat C_n^z$, let
$
\widehat{\sigma}_{n}^z
$
be the total height of all RACS in this clump:
$$
\widehat{\sigma}_{n}^z = \sum_{x\in {\widehat C_n^z}} \sum_{j=1}^{M_{n}^x} \sigma_{n,j}^x
= \sum_{x\in {\widehat C_n^z}} \sigma_n^{x,sum}.
$$
The convention is again that the last quantity is 0 if $\widehat C_n^z=\emptyset$.
We conclude also that $\widehat{\sigma}_{n}^z$ has a light-tailed distribution.

By using monotonicity property 3 (applied to the profile case), one gets that Boolean Model 5, which satisfies the equation
\begin{equation}
\label{eq:5}
\widehat{H}^z_n = \widehat{\sigma}_{n}^z + \max_{y\in \widehat C_n^z \bigcup \{z\}} \widehat{H}^{y}_{n-1},
\end{equation}
with the initial condition $\widehat{H}_0^z=0$ a.s.,
forms an upper bound to Model 4. Similarly,
\begin{equation}
\label{eq:5b}
\widehat{\HH}^z_n = \widehat{\sigma}_{-1}^z +
\max_{y\in \widehat C_{-1}^z \bigcup \{z\}} \widehat{\HH}^{y}_{n-1}\circ \theta^{-1},
\end{equation}
where $\theta$ is the discrete shift on the sequences $\{ \widehat{\sigma}_k^z, \widehat C_{k}^z\}$.
By combining all the bounds constructed so far, we get:
\begin{equation} \label{eq:bornfond}
H^x_t \le \widehat{H}^{z(x)}_{\lceil t\rceil} \quad {\rm and}\quad
\HH^x_t \le \widehat{\HH}^{z(x)}_{\lceil t\rceil} \quad a.s.
\end{equation}
for all $x$ and $t$.

The drawbacks of (\ref{eq:5}) are twofold:
\begin{itemize}
\item[(i)] for all fixed $n$, the random variables $\{\widehat C_n^z\}_z$ are dependent.
This is a major difficulty which will be taken care of by building
a branching upper-bound in subsections \ref{secindset} and \ref{sec:bub} below.
\item[(ii)] for all given $n$ and $z$, the random variables $\widehat C_n^z$
and $\widehat{\sigma}_{n}^z$ are dependent.
We will take care of this by building a second upper bound model in subsection \ref{secinhe} below.
\end{itemize}
Each model will bound (\ref{eq:5}) from above and will hence provide an upper bound
to the initial continuous time, continuous space Poisson hail model.

\subsection{The Branching Upper-bounds}

\subsubsection{The Independent Set Version}\label{secindset}
Assume that the Boolean Model 5 (considered above) has no infinite clump.
Let again $\widehat C_n^x$ be the clump containing
$x\in \Z^d$ at time $n$. For $x\ne y\in \Z^d$, either
$\widehat C_n^x=\widehat C_n^y$ or these two (random) sets are disjoint, which
shows that these two sets are not independent.\footnote{Here
``independence of sets''
has the probabilistic
meaning: two random sets $V_1$ and $V_2$ are {\it independent} if
${\mathbf P} (V_1=A_1,V_2=A_2) = {\mathbf P} (V_1=A_1) {\mathbf P}(V_2=A_2)$, for
all $A_1,A_2\subseteq \Z^d$.} The aim of the following construction
is to show that a certain independent version of these two sets is "larger"
 (in a sense
to be made precise below) than their dependent version.

Below, we call $(\Omega,{\cal F}, {\mathbf P})$ the probability space
that carries the i.i.d. variables
$$\{(\sigma^{z,sum}_0,R^{z,max}_0)\}_{z\in \Z^d},$$
from which the random variables $\{(\widehat C_0^z,\widehat \sigma_0^z)\}_{z\in \Z^d}$ are built.

\begin{lem}\label{lemma001} Assume that
$\lambda < \lambda_c$.
Let $x\ne y$ be two points in $\Z^d$. There exists an extension
of the probability space $(\Omega,{\cal F}, {\mathbf P})$,
denoted by $(\underline \Omega,\underline {\cal F}, \underline {\mathbf P})$, which carries
another i.i.d. family $$\{({\underline \sigma}^{z,sum}_0, {\underline R}^{z,max}_0)\}_{z\in \Z^d}$$
and a random pair $(\widehat{\underline{C}}_0^y,\widehat{\underline{\sigma}}_0^y)$ built from
the latter in the same way as the random variables $\{(\widehat C_0^z,\widehat \sigma_0^z)\}_{z\in \Z^d}$ are built
from $\{(\sigma^{z,sum}_0,R^{z,max}_0)\}_{z\in \Z^d}$, and such that the following properties hold:
\begin{enumerate}
\item The inclusion
$$
\widehat C_0^x\cup \widehat C_0^y \subseteq \widehat C_0^{x} \cup {\widehat {\underline C}}_0^{y},
$$
holds a.s.
\item The random pairs $({\widehat C}_0^{x}, {\widehat \sigma}_0^x)$ and
$({\widehat {\underline C}}_0^{y}, {\widehat {\underline \sigma}}_0^y)$ are independent,
i.e.
\begin{eqnarray*}
& & \hspace{-3cm}\underline {\mathbf P} ({\widehat C}_0^{x}=A_1, {\widehat \sigma}_0^x \in B_1,
{\widehat {\underline C}}_0^{y}=A_2, {\widehat {\underline \sigma}}_0^y \in B_2 )\\
&=&
\underline {\mathbf P}
(\widehat C_0^{x}=A_1, {\widehat \sigma}_0^x \in B_1)
\underline {\mathbf P} ({\widehat {\underline C}}_0^{y}=A_2, {\widehat {\underline \sigma}}_0^y \in B_2)
\\
&=&
{\mathbf P} (\widehat C_0^{x}=A_1, {\widehat \sigma}_0^x \in B_1)
\underline {\mathbf P} ({\widehat {\underline C}}_0^{y}=A_2, {\widehat {\underline \sigma}}_0^y \in B_2),
\end{eqnarray*}
for all sets $A_1,B_1$ and $A_2,B_2$.
\item The pairs $({\widehat {\underline C}}_0^{y}, {\widehat {\underline \sigma}}_0^y)$
and $({\widehat C}_0^{y}, {\widehat \sigma}_0^y)$ have the same law, i.e.
$$ \underline {\mathbf P} ({\widehat {\underline C}}_0^{y}=A,
{\widehat {\underline \sigma}}_0^y \in B)
)= {\mathbf P} ({\widehat {C}}_0^{y}=A, {\widehat \sigma}_0^y \in B),$$
for all sets $(A,B)$.
\end{enumerate}
\end{lem}

\paragraph{Proof.}

We write for short $\widehat{C}^x = \widehat{C}^x_0$ and $\widehat{\sigma}^x =
\widehat{\sigma}^x_0$.
Consider first the case of balls with a constant integer radius $R=R^{max}$ (the case with
random radii is considered after).
Recall that we consider $L_\infty$-norm balls
in ${\R}^d$, i.e. $d$-dimensional cubes with side $2R$, so
a "ball $B^x$ centered at point $x=(x_1,\ldots , x_d)$" is the closed cube $x+ [-R,+R]^d$.

We assume that the ball $B^x$
exists at time 0 with probability $p = {\mathbf P}(M\ge 1)\in (0,1)$ independently
of all the others. Let $E^x=B^x$ if $B^x$ exists at time 0 and $E^x=\emptyset$, otherwise,
and let $\alpha^x = {\bf I} (E^x=B^x)$ be the indicator of the event that $B^x$ exists
(we drop the time index to have lighter notation).
Then the family of r.v.'s $\{\alpha^x\}_{x\in \Z^d}$ is i.i.d.

Recall that the clump $\widehat C^x$, for the input $\{\alpha^x\}$, is the maximal connected set of
balls that contains $x$. This clump is empty if and only if $\alpha^y=0$, for all $y$ with
$d_{\infty}(x,y)\le R$. Let $L(\widehat C^x)$ denote the number of lattice points in the clump $\widehat C^x$,
$0\le L(\widehat C^x) \le \infty$. Clearly, $L(\widehat C^x)$ forms a stationary (translation-invariant) sequence.

For all sets $A \subset {\Z}^d$, let
$$Int(A) = \{ x\in A : \ B^x \subseteq A\} ,\quad \mbox{and}\quad
Hit(A)= \{ x \in \Z^d : \ B^x \cap A \ne \emptyset \}.$$
For $A$ and $x,y\in A$, we say that the event
$$
\Bigl\{ x\ { \Longleftrightarrow  \atop {Int(A), \{\alpha^u\}}}\  y\Bigr\}
$$
occurs if, for the input $\{\alpha^u\}$, the random set $E^A = \bigcup_{z\in Int(A)} E^z$ is
connected and both $x$ and $y$ belong to $E^A$.

Then the following events are equal:
$$
\Bigl\{\widehat C^{x} =A\Bigr\} = \bigcap_{z\in A}
\Bigl\{ x\ { \Longleftrightarrow\atop {Int(A), \{\alpha^u\}}}\ z\Bigr\}
\bigcap \bigcap_{z\in Hit(A) \setminus Int(A)} \{\alpha^{z}=0\}.
$$
Therefore, the event $\{\widehat C^{x}=A \}$ belongs to the sigma-algebra ${\cal F}^{\alpha}_{Hit (A)}$
generated by the random variables $\{\alpha^{x}, x \in Hit (A)\}$. Let also
 ${\cal F}^{\alpha ,\sigma}_{Hit (A)}$ be the sigma-algebra generated by the random
variables  $\{\alpha^{x}, \sigma^x,  x \in Hit (A)\}$.

Recall the notation $\sigma^{z,sum}_0 = \sum_{j=1}^{M_{0}^z}\sigma_{0,j}^z$.
We will write for short $\sigma^z = \sigma_0^{z,sum}$. Clearly $\sigma^z=0$ if $\alpha^z=0$, and
the family of pairs $\{(\alpha^z, \sigma^z)\}$ is i.i.d. in $z\in \Z^d$.

Let $\{(\alpha^z_{*}, \sigma^z_{*})\}$ be another i.i.d. family in $z\in \Z^d$ which does not depend
on all random variables introduced earlier and whole elements have a common distribution
with $(\alpha^0, \sigma^0)$.
Let  $(\underline \Omega,\underline {\cal F}, \underline {\mathbf P})$
be the product probability space that carries both  $\{(\alpha^{z},\sigma^z)\}$ and
$\{(\alpha^{z}_*,\sigma^z_*)\}$.
Introduce then a third family $\{({\underline{\alpha}}^{z},{\underline{\sigma}}^z)\}$ defined as
follows: for any set $A$ containing $x$, on the event $\{\widehat C^{x}=A\}$ we let
\begin{eqnarray*}
({\underline{\alpha}}^{z}(A), {\underline{\sigma}}^z(A)) =
\begin{cases}
 (\alpha^{z}_{*},\sigma^z_{*}) & \mbox{if} \quad z\in Hit (A) \\
(\alpha^{z},\sigma^z), & \mbox{otherwise}.
\end{cases}
\end{eqnarray*}
When there is no ambiguity, we will use the notation
$({\underline{\alpha}}^{z}, {\underline{\sigma}}^z)$
in place of 
$({\underline{\alpha}}^{z}(A), {\underline{\sigma}}^z(A))$.
First, we show that $\{({\underline {\alpha}}^z, {\underline{\sigma}}^z)\}$ is an i.i.d. family.
Indeed, for any finite
set of distinct points $y_1,\ldots,y_k$, for any $0-1$-valued sequence $i_1,\ldots,i_k$,
and for all measurable sets $B_1,\ldots ,B_k$,
\begin{eqnarray*}
\underline {\mathbf P} (\underline{\alpha}^{y_j} = i_j, \underline{\sigma}^{y_j}\in B_j, j=1,\ldots,k)
&=&
\sum_A \underline {\mathbf P}
(\widehat C^{x}=A,  \underline{\alpha}^{y_j} = i_j, \underline{\sigma}^{y_j}\in B_j, j=1,\ldots,k)\\
& & \hspace{-6.5cm}=
\sum_A \underline {\mathbf P}
(\widehat C^{x}=A, {\alpha}^{y_j}_* = i_j, \sigma^{y_j}_*\in B_j, y_j\in Hit (A) \
\mbox{and}
\ \alpha^{y_j}=i_j, \sigma^{y_j}\in B_j, y_j \in (Hit (A))^c)\\
& & \hspace{-6.5cm}=
\sum_A \underline {\mathbf P}
(\widehat C^{x}=A) \underline {\mathbf P}( {\alpha}^{y_j}_* = i_j, \sigma^{y_j}_*\in B_j, y_j\in Hit (A) )
\underline {\mathbf P} (\alpha^{y_j}=i_j, \sigma^{y_j}\in B_j, y_j \in (Hit (A))^c)\\
& & \hspace{-6.5cm}= \sum_A {\mathbf P} (\widehat C^{x}=A)
 \prod_{j=1}^k {\mathbf P} (\alpha^0=i_j, \sigma^0\in B_j)
\\
& & \hspace{-6.5cm}=
 \prod_{j=1}^k {\mathbf P} (\alpha^0=i_j, \sigma^0\in B_j).
\end{eqnarray*}
Notice that the sum over $A$ is a sum over finite $A$. This keeps the number of terms
countable. This is licit due to assumption on the finiteness of the clumps.

Let $\widehat {\underline{C}}^{y}$ be the clump of $y$ for $\{\underline {\alpha}^z\}$ and let
$\widehat {\underline \sigma}^y = \sum_{z\in  \widehat {\underline{C}}^{y}}
{\underline{\sigma}}^z$.
We now show that the pairs $({\widehat C}^{x}, {\widehat \sigma}^x)$ and
$({\widehat {\underline{C}}}^{y}, {\widehat {\underline \sigma}}^y)$ are independent.
For all sets $A$, let ${\cal F}^A$ be the sigma-algebra
generated by the random variables
$$(\alpha^{(A)}, \sigma^{(A)}) =
\{(\alpha^u_{*},\sigma^u_*, u\in Hit (A); \alpha^v, \sigma^v, v\in (Hit (A))^c\},$$
and let $\widehat {\underline C}^{y}(A)$
be the clump containing $y$
in the environment $\alpha^A$. Let also
 $  \widehat {\underline \sigma}^y (A) = \sum_{z\in  \widehat {\underline{C}}^{y}}
{\underline{\sigma}}^z(A)$. Clearly, $(\alpha^{(A)},\sigma^{(A)})$ is also an i.i.d. family.
Then, for all sets $A_1,B_1$ and $A_2,B_2$,
\begin{eqnarray}
& & \hspace{-3cm}  \underline {\mathbf P} (\widehat C^{x}=A_1, \widehat \sigma^x\in B_1,
\widehat {\underline C}^{y} = A_2, \widehat {\underline \sigma}^y\in B_2)\nonumber\\
& =&
\underline {\mathbf P} (\widehat C^{x}=A_1,  \widehat \sigma^x\in B_1,
\widehat {\underline C}^{y}(A_1)=A_2, \widehat {\underline \sigma}^y(A_1) \in B_2) \nonumber \\
&=&
\underline {\mathbf P} (\widehat C^{x}=A_1, \widehat \sigma^x\in B_1)
\underline  {\mathbf P} (\widehat C^{y}(A_1)=A_2, \widehat \sigma^y (A_1) \in B_2)\nonumber \\
&=&
{\mathbf P} (\widehat C^{x}=A_1, \widehat \sigma^x\in B_1)
{\mathbf P} (\widehat{C}^{y}=A_2, \widehat \sigma^y\in B_2).
\end{eqnarray}
The second equality follows from the fact that
the event $\{\widehat C^{x}=A_1, \widehat \sigma^x\in B_1\}$ belongs to the sigma-algebra ${\cal F}^{\alpha , \sigma}_{Hit (A_1)}$
whereas the event $\{\widehat {\underline C}^{y}(A_1)=A_2,
\widehat {\underline \sigma}^y(A_1) \in B_2\}$ belongs to
the sigma-algebra ${\cal F}^{A_1}$, which is independent.
The last equality follows from the fact that $\{\alpha^{(A_1)}, \sigma^{(A_1)}\}$ is
an i.i.d. family with
the same law as $\{\alpha^x, \sigma^x\}$.

We now prove the first assertion of the lemma.
If $\widehat C^{x} =\widehat C^{y}$,
then the inclusion is obvious. Otherwise,
$\widehat C^{x}\bigcap \widehat C^{y} = \emptyset$ and if $\widehat C^{x}=A$,
the size and the shape of $\widehat C^{y}$ depend only on
$\{\alpha^u, u\in (Hit (A))^c\}$. Indeed, on these events,
$$v\in \widehat C^{y}\quad {\rm iff}\quad
y \ {{\Longleftrightarrow}\atop{Int (A^c), \{\alpha^x\}}}\  v\  .$$
Then the first assertion follows since, first,
the latter relation is determined by
$\{\alpha^u, u\in Int (A^c) \}$ and, second, $Int (A^c) = (Hit (A))^c$. We
may conclude that
$\widehat {\underline C}^{y}(A) \supseteq \widehat C^{y}$ because
some $\alpha_*^{z}, z\in Hit (A) \setminus Int (A)$ may take value 1.

Finally, the second assertion of the lemma follows from the construction.

The proof of the deterministic radius case is complete.

Now we turn to the proof in the case of random radii. Recall that we assume that
the radius $R$ of a Model 2 RACS is a positive integer-valued r.v. and this is a radius in the
$L_{\infty}$ norm. For $x\in {\Z}^d$ and $k=1,2,\ldots$,
let $B^{x,k}$ be  the $L_{\infty}$-norm ball with center $r$ and radius $k$. Recall that
$M_0^{x,k}$ is the number of RACS that arrive at time $0$, are centered at $x$ and
have radius $k$. Then, in particular,
$$
R^{x,max}_0 = \max \{ k : \ M_0^{x,k}\ge 1\}.
$$

Let $\alpha^{x,k}$ be the indicator of event $\{ M_0^{x,k}\ge 1\}$ and $E^{x,k}$ a random set,
$$
E^{x,k}=B^{x,k} \quad \mbox{if} \quad \alpha^{x,k}=1 \quad \mbox{and}
\quad E^{x,k}=\emptyset, \quad \mbox{otherwise}.
$$
Again, the r.v.'s
$\alpha^{x,k}$ are mutually independent (now both in $x$
and in $k$) and also i.i.d. (in $x$).

For each $A\subseteq {\Z}^d$, we let $Int_2 (A) = \{ (x,k): \ x\in A, k\in {\N},
B^{x,k}\subseteq A\}$ and $Hit_2 (A) = \{ (x,k): \ x\in A, k\in {\N}, B^{x,k} \bigcap A
\ne \emptyset \}$.

For $x,y\in A$, we say that the event
$$
\Bigl\{ x\ {\Longleftrightarrow\atop {Int_2(A), \{\alpha^{u,l}\}}}\  y\Bigr\}
$$
occurs if, for the input $\{\alpha^x\}$, the random set $E^A = \bigcup_{(z,k)\in Int(A)} E^{z,k}$ is
connected and both $x$ and $y$ belong to $E^A$.

Then the following events are equal:
$$
\Bigl\{\widehat C^{x} =A\Bigr\} = \bigcap_{z\in A}
\Bigl\{ x\ {\Longleftrightarrow\atop {Int_2(A), \{\alpha^{u,l}\}}}\ z\Bigr\}
\bigcap \bigcap_{(z,k)\in Hit_2(A) \setminus Int_2(A)} \{\alpha^{z,k}=0\}.
$$
Therefore, the event $\{\widehat C^{x}=A \}$ belongs to the sigma-algebra ${\cal F}^{\alpha}_{Hit_2 (A)}$
generated by the random variables $\{\alpha^{x,k}, (x,k) \in Hit_2 (A)\}$.
For $x\in \Z^d$ and for $k=1,2,\ldots$, we let $\sigma^{x,k}=\sum_{j=1}^{M_0^{x,k}} \sigma^x_{0,j}$
where the sum of the heights is taken over all RACS that arrive at time $0$, are centered at $x$ and
have radius $k$. Clearly, the random vectors $(\alpha^{n,k},\sigma^{n,k})$ are independent in all
$x$ and $k$ and identically distributed in $x$, for each fixed $k$.

Let $\{(\alpha^{x,k}_*,\sigma^{x,k}_*)\}$ be another independent family of pairs
that does not depend on all random
variables introduced earlier and is such that, for each $k$ and $x$, the pairs
$(\alpha^{x,k}_*,\sigma^{x,k}_*)$ and
$(\alpha^{0,k},\sigma^{0,k})$ have a common distribution.
Let $(\underline \Omega,\underline {\cal F}, \underline {\mathbf P})$
be the product probability space that carries both  $\{(\alpha^{x,k},\sigma^{x,k})\}$ and
$\{(\alpha^{x,k}_*,\sigma^{x,k}_*)\}$.
Introduce then a third family $\{(\underline{\alpha}^{x,k},\underline{\sigma}^{x,k})\}$ defined as
follows: for any set $A$ containing $x$, on the event $\{\widehat C^{x}=A\}$ we let
\begin{eqnarray*}
(\underline{\alpha}^{z,l},\underline{\sigma}^{z,l}) =
\begin{cases}
 (\alpha^{z,l}_{*},\sigma^{z,l}_*) & \mbox{if} \quad (z,l)\in Hit_2 (A) \\
(\alpha^{z,l},\sigma^{z,l}), & \mbox{otherwise}.
\end{cases}
\end{eqnarray*}
The rest of the proof is then quite similar to that of the constant radius case:
we introduce again $\widehat {\underline{C}}^{y}$, which is now
the clump of $y$ for $\{\underline {\alpha}^{z,l}\}$ with the height
$\widehat {\underline \sigma}^y = \sum_k \sum_{z\in \widehat C^y}
{\underline \sigma}^{z,k}$; we
then show that the random pairs  $(\widehat C^{x},\widehat \sigma^x)$ and
$(\widehat {\underline{C}}^{y}, \widehat {\underline \sigma}^y)$ are independent
and finally establish the first and the second assertions of the lemma.

\hfill $\Box$

We will need the following two remarks on Lemma \ref{lemma001}.

\begin{rem} \label{rem2}
In the proof of Lemma \ref{lemma001}, the roles of the points $x$ and $y$ and of the sets $\widehat{C}^x$
and $\widehat{C}^y$ are not symmetrical. It is important that $\widehat{C}^x$ is a clump while
from $V= \widehat{C}^y$, we only need the following monotonicity property:
the set $V\setminus \widehat{C}^x$ is a.s. bigger in the environment
$\{\underline{\alpha}^{z}\}$ than in the environment $\{\alpha^z\}$. One can note that
any finite union of clumps also satisfies this last property.
\end{rem}

\begin{rem} \label{rem22}
From the proof of Lemma \ref{lemma001}, the following properties hold.
\begin{enumerate}
\item
On the event where $\widehat C_0^x$ and  $\widehat C_0^y$ are
disjoint, we have  $\widehat{C}_0^y \subseteq \widehat{\underline C}_0^y$
and $\sigma^{z,sum}_0 = {\underline \sigma}^{z,sum}_0$ a.s.,
for all $z\in \widehat C_0^y$, so that $\widehat{\sigma}_0^y\le \widehat{\underline \sigma}_0^y$.
\item
On the event where
$\widehat{C}_0^x=\widehat{C}_0^y$, we have
$\widehat{\sigma}_0^x=\widehat{\sigma}_0^y$.
\end{enumerate}
Let us deduce from this that, for all constants $a^x\ge a^y$,
for all $z\in \widehat{C}_0^{x}\cup \widehat{C}_0^{y}$,
there exists a random variable $r(z) \in \{x,y\}$ such that
$z\in \widehat{\underline C}_0^{r(z)}$ (with the convention
$\widehat{\underline C}_0^{x}=\widehat{C}_0^{x}$ and
$\widehat{\underline \sigma}_0^{x}=\widehat{\sigma}_0^{x}$) a.s. and
$$\max_{u\in\{x,y\} \ :\ z \in \widehat{C}^{u}_0}
\left( a^{u} + \widehat{\sigma}_0^{u}\right)
\le a^{r(z)} + \widehat{\underline \sigma}_0^{r(z)} \quad \mbox{a.s.}$$
In case 2 and case 1 with $z\in \widehat{C}_0^{x}$,
we take $r(z)=x$ and use the fact that $a^x\ge a^y$.
In case 1 with $z\in \widehat{C}_0^{y}$, we take $r(z)=y$ and use the fact that
$\widehat{\sigma}_0^y\le \widehat{\underline \sigma}_0^y$.

As a direct corollary of the last property, the inequality
$$
\max (a^x+{\widehat \sigma}_0^x, a^y+{\widehat \sigma}_0^y)
\le
\max (a^x+{\widehat \sigma}_0^x, a^y+{\widehat {\underline \sigma}}_0^y)
$$
holds a.s.
Here $\widehat{\underline \sigma}_0^y =
\sum_{z\in \widehat{\underline C}_0^y}
\underline{\sigma}_0^{z,sum}$.
\end{rem}

We are now in a position to formulate a more general result:

\begin{lem}\label{general0}
Assume again that $\lambda < \lambda_c$. 
Let $S$ be a set of $\Z^d$ of cardinality $p\ge 2$, say $S=\{x_1,\ldots ,x_{p}\}$.
There exists an extension of the initial probability space
and random pairs ${(\widehat {\underline C}}_0^{x_i},\widehat {\underline \sigma}_0^{x_i})$,
$i=2,\ldots,p$
defined on this extension which are such that:
\begin{enumerate}
\item The inclusion
\begin{equation}\label{eq:lem2p1}
\bigcup_{j=1}^{p} \widehat C_0^{x_j} \subseteq
\bigcup_{j=1}^{p}
\widehat {\underline  C}^{x_{j}}_0 \quad
a.s.
\end{equation}
holds with $\widehat C_0^{x_1} =\widehat {\underline  C}^{x_{1}}_0$.
\item
For all real valued constants $a^{x_1},a^{x_2}, \ldots , a^{x_{p}}$ such
that $a^{x_1}= \max_{1\le i \le p} a^{x_i}$,
for all $z\in \bigcup_{j=1}^p \widehat{C}_0^{x_j}$,
there exists a random variable $r(z) \in \{x_1,\ldots,x_p\}$ such that
$z\in \widehat{\underline C}_0^{r(z)}$ a.s. and
\begin{equation}\label{eq:lem2p21}
\max_{j\in \{1,\ldots,p\} \ : \ z \in \widehat{C}^{x_j}_0}
\left( a^{x_j} + \widehat{\sigma}_0^{x_j}\right)
\le a^{r(z)} + \widehat{\underline \sigma}_0^{r(z)} \quad \mbox{a.s.}
\end{equation}
In particular, the inequality
\begin{equation}\label{eq:lem2p22}
\max_{1\le j\le p}\left( a^{x_j}+ \widehat{\sigma}^{x_j}_0 \right) \le
\max_{1\le j \le p}
\left(a^{x_j}+\widehat{\underline \sigma}^{x_j}_0
\right)
\end{equation}
holds a.s. with $\widehat \sigma_0^{x_1} =\widehat {\underline  \sigma}^{x_{1}}_0$.
\item The pairs $ (\widehat C_0^{x_1},\widehat \sigma^{x_1}_0),
(\widehat {\underline  C}^{x_{2}}_0,\widehat {\underline \sigma}_0^{x_{2}}),
\ldots,
(\widehat {\underline  C}^{x_{p}}_0,\widehat {\underline \sigma}_0^{x_{p}})$
are mutually independent.
\item The pairs $( \widehat C_0^{x_i},\widehat \sigma^{x_i}_0)$ and
$(\widehat {\underline C}_0^{x_i},\widehat {\underline \sigma}^{x_i}_0)$, have
the same law, for each fixed $i=2,\ldots,p$. 
\end{enumerate}
\end{lem}

\paragraph{Proof.}
We proceed by induction on $p$. 
Assume the result holds for any set with $p$ points. Then consider a
set $S$ of cardinality $(p+1)$ and number its points arbitrarily,
$S=\{x_1,\ldots,x_{p+1})$.
For $A$ fixed, consider the event
$\{\widehat{C}^{x_1}_0=A\}$. On this event,
define the same family $(\underline \alpha^{z,l},
\underline \sigma^{z,l})$ as in the previous proof and
consider the $p$ clumps $\underline D^{x_2},\ldots,\underline D^{x_{p+1}}$ with their heights, say
$\underline{s}^{x_2},\ldots , \underline{s}^{x_{p+1}}$ for this family.
By the same reasons as in the proof of Lemma \ref{lemma001}, $(\widehat{C}^{x_1}_0,\sigma^{x_1}_0)$ is independent of
$(\underline D^{x_2},\underline s^{x_2}),\ldots,(\underline D^{x_{p+1}},\underline{s}^{x_{p+1}})$.
By Remark \ref{rem2},
$$
\bigcup_{j=1}^{p+1} \widehat C_0^{x_j} \subseteq
\widehat{C}_0^{x_1} \bigcup
\bigcup_{j=2}^{p+1} \underline D^{x_{j}}\quad
{a.s.}
$$
By the induction step,
$$ \underline D^{x_2}\cup \ldots\cup \underline D^{x_{p+1}}\subseteq_{a.s.} \widehat {\underline  C}^{x_{2}}_0\cup
\ldots \cup \widehat {\underline  C}^{x_{p+1}}_0,$$
with $\widehat {\underline  C}^{x_{2}}_0, \ldots, \widehat {\underline  C}^{x_{p+1}}_0$ defined
as in the lemma's statement
and then the first, third and fourth assertions follow.

We now prove the second assertion, again by induction on $p$. If $p=2$, this is Remark \ref{rem22}.
For $p>2$, we define $L_1 = \{ p+1\ge j\ge 1 \ : \widehat{C}_0^{x_j}=\widehat{C}_0^{x_1} \}$
and we consider two cases:
\begin{enumerate}
\item $z\in \widehat { C}^{x_{1}}_0$. In this case
let $\overline{L}_1=\{1, \ldots , p+1, \}\setminus L_1$.
Since $z\notin \widehat C_0^{x_j}$ for $j\in \overline{L}_1$ and since $\widehat \sigma_0^{x_j}=\widehat \sigma_0^{x_1}$ for all
$j\in L_1$, we get that (\ref{eq:lem2p21}) holds with $r(z)=1$ when using the fact that $a^{x_1}= \max_{1\le i \le p} a^{x_i}$.
\item $z\notin \widehat { C}^{x_{1}}_0$. In this case
let $\overline{L}_1^z=\{ 1\le j\le p+1\ : \ j\notin L_1, z\in \widehat{C}^{x_j}\}$. We can assume w.l.g. that this set is non-empty.
Then for all $j\in \overline{L}_1^z$, we have  $\underline s^{x_j}\ge\widehat{\sigma}^{x_j}$, by Lemma \ref{lemma001} and
Remark \ref{rem2}.  So
$$
\max_{j\in \overline{L}_1^z} \left(a^{x_j}+\widehat{\sigma}_0^{x_j} \right) \le
\max_{j\in\overline{L}_1^z} \left(a^{x_j}+\underline s^{x_j} \right)
\quad \mbox{a.s.}
$$
Now, since the cardinality of $\overline L_1^z$ is less than or equal to $p$, we can use the induction assumption,
which shows that when choosing $i_1\in\overline{L}_1^z$
such that $a^{x_{i_1}} = \max_{i\in\overline{L}_1^z}a^{x_i}$,
we have
$$ \max_{j\in\overline{L}_1^z} \left(a^{x_j}+\underline s^{x_j} \right)
\le a^{x_{r(z)}}+
\widehat{ \underline \sigma} _0^{x_{r(z)}}
,$$
with $r(z)\in \overline{L}_1^z$ and
with the random variables $\{\widehat{\underline \sigma}_0^{x_j}\}$ defined as in the lemma's statement,
but for $\widehat{\underline \sigma}_0^{x_{i_1}}$ which we take equal to $\underline s^{x_{i_1}}$.
The proof in concluded in this case too when using the fact that
the random variable $\underline s^{x_{i_1}}$ is mutually independent of the random variables
$(\{\widehat{ \underline \sigma} _0^{x_j}\},\widehat{\sigma}_0^{x_1})$ and it has the
same law as $\widehat{\underline \sigma}^{x_{i_1}}$.
\end{enumerate}
\hfill $\Box$

\subsubsection{Comparison with a Branching Process}
\label{sec:bub}

\paragraph{Paths and Heights in Boolean Model 5}
Below, we focus on the backward construction associated with Boolean Model 5,
for which we will need more notation.

Let $\D_n^x$ denote the set of {\em descendants} of level $n$ of
$x\in \R^d$ in this backward process, defined as follows:
\begin{eqnarray*}
\D_1^x & = & \widehat C_0^x \cup \{x\} \\
\D_{n+1}^x & = & \bigcup\limits_{y\in \D_{n}^x} \widehat C_{-n}^y \cup \{y\}, \quad n\ge 1.
\end{eqnarray*}
By construction, $\D_n^x$ is a non-empty set for all $x$ and $n$.
Let $d_n^x$ denote the cardinality of  $\D_n^x$.

Let $\Pi_n^x$ denote the set of paths starting from $x=x_0\in \Z^d$ and
of length $n$ in this backward process: $x_0,x_1,\ldots,x_n$ is
such a path if $x_0,x_1,\ldots,x_{n-1}$ is a path of length $n-1$
and $x_n\in \widehat C_{-n+1}^{x_{n-1}} \cup \{x_{n-1}\}$.
Let $\pi_n^x$ denote the cardinality of  $\Pi_n^x$. Clearly,
$d_n^x\le\pi_n^x$ a.s., for all $n$ and $x$.

Further, the {\it height} of a path $l_n = (x_0,\ldots,x_n)$ is the sum of the heights
of all clumps along the path:
$$
\sum_{i=0}^{n-1} {\widehat \sigma}_{-i}^{x_i}.
$$
In particular, if the paths $l_n$ and $l_n^{'}$ differ only by the last points
$x_n\in \widehat{\sigma}_{-n+1}^{x_{n-1}}$ and $x'_n\in \widehat{\sigma}_{-n+1}^{x_{n-1}}$,
then their heights coincide.

For $z\in \Z^d$, let $\widehat h_n^{x,z}$ be the maximal height of all paths of length $n$ that start from
$x$ and end at $z$, where the maximum over the empty set is zero.

Let $\widehat{\HH}_n^x$, $n\ge 0$ be the maximal height of all paths of length $n$ that start from
$x$.
Then $\widehat{\HH}(n)=\max_z \widehat h_n^{x,z}$.

\paragraph{Paths and Heights in a Branching Process}
Now we introduce a branching process (also in the backward time) that starts from point $x=x_0$ at
generation 0. Let $(V_{n,i}^z,s_{n,i}^{z})$, $z\in \Z^d$, $n\ge 0$, $i\ge 1$ be a family of mutually
independent random pairs such that, for each $z$, the pair
$(V_{n,i}^z,s_{n,i}^{z})$ has the same distribution
as the pair $(\widehat{C}_0^z\cup\{z\},\widehat \sigma_0^z)$, for all $n$ and $i$.

In the branching process defined below, we do not distinguish between points and paths.

In generation 0, the branching process has one point: $\widetilde{\Pi}_0^{x_0}= \{(x_0\}$.
In generation 1, the points of the branching process are
$\widetilde{\Pi}_1^{x_0}= \{(x_0,x_1), \  x_1\in V_{0,1}^{x_0}\}$.
Here the cardinality of this set is the number of
points in $V_{0,1}^{x_0}$ and all end coordinates $x_1$ differ (but this is not the case
for $n\ge 2$, in general).

In generation 2, the points of the branching process are
$$
\widetilde{\Pi}_2^{x_0}= \{(x_0,x_1,x_2), \ (x_0,x_1)\in \widetilde{\Pi}_1^{x_0},
x_2\in V_{1,1}^{x_1}\}.
$$
Here a last coordinate $x_2$ may appear several times, so we introduce a
multiplicity function $k_2$: for $z\in\Z^d$,
$k_2^z$ is the number of $(x_0,x_1,x_2)\in\widetilde{\Pi}_1^{x_0}$ such that
$x_2=z$.

Assume the set of all points in generation $n$ is $\widetilde{\Pi}_n^{x_0}=\{ (x_0, x_1,\ldots,x_n)\}$
and $k_n^z$ is the multiplicity function (for the last coordinate).
For each $z$ with $k_n^z >0$, number arbitrarily
all points with last coordinate $z$ from 1 to $k_n^z$
and let $q(x_1,x_2,\ldots,x_n)$ denote the number
given to point $(x_0,\ldots,x_n)$ with $x_n=z$.
Then the set of points in generation $n+1$ is
$$
\widetilde{\Pi}_{n+1}^{x_0}=
\{(x_0,\ldots,x_n,x_{n+1}\ : \
(x_0,\ldots,x_n)\in \widetilde{\Pi}_n^{x_0}, \ x_{n+1}\in V_{n,q(x_0,\ldots,x_n)}^{x_n} \}.
$$
Finally the height of point
$(x_0,\ldots,x_n)\in\widetilde{\Pi}_n^{x_0}$ is defined as
$$
\widetilde{h}(x_0,\ldots,x_n)= \sum_{i=0}^{n-1} s_{i,q_i}^{x_i}~,
$$
where $q_i=q(x_0,\ldots,x_i)$.

\paragraph{Coupling of the two Processes}
\begin{lem}\label{lemmanew}
Let $x_0$ be fixed. Assume that $\lambda <\lambda_c$. 
There exists a coupling of Boolean Model 5 and of the branching process 
defined above such that,
for all $n$,
for all points $z$ in the set $\D_n^{x_0}$, there exists a point
$(x_0,\ldots,x_n)\in\widetilde{\Pi}_n^{x_0}$ such that $x_n=z$ and
$\widehat h_n^{x_0,z}\le\widetilde{h}(x_0,\ldots,x_n)$ a.s.
\end{lem}

\paragraph{Proof}
We construct the coupling and prove the properties by induction. For $n=0,1$, the process of Boolean Model 5
and the branching process coincide. Assume that the statement of the lemma
holds up to generation $n$. For $z\in\D_n^{x_0}$, let $a^z=\widehat h_n^{x_0,z}$.

Now, conditionally on the values of both processes up to level $n$ inclusive,
we perform the following coupling at level $n+1$:
we choose $z_*$ with the maximal $a^z$ and we apply Lemma \ref{general0} with $S=\D_n^{x_0}$,
with $z_*$ in place of $x_1$, and with $\{\widehat{C}_{-n}^z\}_z$
(resp. $\{\widehat{\underline C}_{-n}^z\}_z$) in place of   $\{\widehat{C}_{0}^z\}_z$
(resp. $\{\widehat{\underline C}_{0}^z\}_z$); we then take
\begin{itemize}
\item $V_{n,1}^{z_*}=\widehat{C}_{-n}^{z_*}\cup\{{z_*}\}$;
\item $V_{n,1}^z=\widehat{\underline C}_{-n}^z\cup\{z\}$ for all $z\in \D_n^{x_0}$, $z\ne z_*$;
\item $s_{n,1}^{z_*}=\widehat{\sigma}_{-n}^{z_*}$;
\item $s_{n,1}^z=\widehat{\underline \sigma}_{-n}^z$ for all $z\in \D_n^{x_0}$, $z\ne z_*$.
\end{itemize}

By induction assumption, for all $z\in \D_n^{x_0}$, there exists a
$(x_0,\ldots,x_n)\in\widetilde{\Pi}_n^{x_0}$ such that $x_n=z$.
This and Assertion 1 in Lemma \ref{general0} show that if
$u\in \D_{n+1}^{x_0}$, then $(x_0,\ldots,x_n,u)\in\widetilde{\Pi}_{n+1}^{x_0}$,
which proves the first property.

By a direct dynamic programming argument, for all $u\in D_{n+1}^{x_0}$
$$
\widehat h_{n+1}^{x_0,u} = \max_{z\in \D_n^{x_0}, u\in \widehat{C}_{-n}^z} \widehat h_{n}^{x_0,z} + \widehat{\sigma}_{-n}^z.
$$
We get from Assertion 2 in Lemma \ref{general0}
applied to the set $\{x_1,\ldots,x_p\}=\{z\in \D_n^{x_0}, u\in \widehat{C}_{-n}^z\}$ that
$$
\widehat h_{n+1}^{x_0,u} \le
\max_{z\in \D_n^{x_0}, u\in \widehat{C}_{-n}^z} \left( \widehat h_{n}^{x_0,z} + \widehat{\underline \sigma}_{-n}^z\right).
$$
By the induction assumption, for all $z$ as above,
$\widehat h_n^{x_0,z}\le\widetilde{h}(x_0,\ldots,x_n)$ a.s. for some $(x_0,\ldots,x_n)\in \widetilde{\Pi}_n^{x_0}$ with $x_n=z$.
Hence for all $u$ as above, there exists a path $(x_0,\ldots,x_n,x_{n+1})\in \widetilde{\Pi}_{n+1}^{x_0}$ with $x_{n+1}=u$
and such that
$\widehat h_{n+1}^{x_0,u}\le \widetilde{h}(x_0,\ldots,x_n,x_{n+1})$ with $(x_0,\ldots,x_n,x_{n+1})\in \widetilde{\Pi}_{n+1}^{x_0}$
and $x_{n+1}=u$.
\hfill $\Box$

\subsubsection{Independent Heights}\label{secinhe}
Below, we assume that the light tail assumptions on $\xi^d$ and $\sigma$ are satisfied
(see Section \ref{ssds}).

In the last branching process, the pairs $(V_{n,i}^z,s_{n,i}^{z})$ are mutually independent in $n,i$ and $z$. However,
for all given $n, i$ and $z$, the random variables $(V_{n,i}^z,s_{n,i}^{z})$,
are dependent. It follows from Proposition \ref{Sergey} in the appendix that
one can find random variables $(W_{n,i}^z,t_{n,i}^{z})$ such that
\begin{itemize}
\item For all $n,i$ and $z$, $V_{n,i}^z\subset W_{n,i}^z$ a.s.
\item The random sets $W_{n,i}^z$ are of the form $z+w_{n,i}^z$, where the sequence $\{w_{n,i}^z\}$
is i.i.d. in $n,i$ and $z$.
\item The random variable ${\rm card} (W_{0,1}^0)$ has exponential moments.
\item For all $n,i$ and $z$, $s_{n,i}^z\le t_{n,i}^z$ a.s.
\item The random variable $t_{0,1}^0$ has exponential moments.
\item The pairs $(W_{n,i}^z,t_{n,i}^{z})$ are mutually independent in $n,i$ and $z$.
\end{itemize}
So the branching process built from the $\{(W_{n,i}^z,t_{n,i}^{z})\}$ variables
is an upper bound to the one built from the $\{(V_{n,i}^z,s_{n,i}^{z})\}$ variables.

\subsection{Upper Bound on the Growth Rate}
The next theorem, which pertains to branching process theory, is not new (see e.g. \cite{Biggins}).
We nevertheless give a proof for self-containedness. It features a
branching process with height (in the literature, one also says with age or with height), starting
from a single individual, as the one defined in Section \ref{secinhe}.
Let $v$ be the typical progeny size, which we assume to be light-tailed.
Let $s$ be the typical height of a node, which we also assume to be light-tailed.
\begin{theorem}\label{LinearGrowth} 
Assume that $\lambda < \lambda_c$. 
For $n\ge 0$,
let $h(n)$ be the maximal height of all descendants of generation $n$ in the branching
process defined above. There exists a finite and positive constant $c$ such that
\begin{equation}
\lim\sup_{n\to\infty} \frac{\widehat\HH(n)} n \le c \quad a.s.
\end{equation}
\end{theorem}

{\sc Proof.}
Let $(v_i,s_i)$ be the i.i.d. copies of $(v,s)$.
Take any positive $a$. Let
$D(a)$ be the event
$$ D(a)=\bigcup\limits_{n\ge 1} \{{d}_n>a^n \},$$
with $d_n$ the number of individuals of generation $n$ in the branching process.
For all $c>0$ and all positive integers $k$,
let $W_{c,k}$ be the event $
\bigl\{ \frac{h(k)} k \le c \bigr\}$.
Then
$$
W_{c,k} \subseteq \bigl(W_{c,k}\cap \overline D(a)\bigr) \bigcup D(a),
$$
where $\overline{D}(a)$ is the complement of $D(a)$.
From Chernoff's inequality, we have, for $\gamma\ge 0$
\begin{eqnarray*}
{\mathbf P}(D(a)) &=&
{\mathbf P} \left(
\bigcup_{n\ge 0}\{ {d}_{n+1} > a^{n+1},
{d}_i \le a^i, \forall i\le n \} \right)\\
&\le &
\sum_{n\ge 1} {\mathbf P}
\left(\sum_{j=1}^{a^n} v_j >a^{n+1}\right) \\
&\le &
\sum_{n\ge 1} \left(
{\mathbf E} \exp (\gamma v ) \right)^{a^n} \cdot e^{-\gamma a^{n+1}}\\
&\le &
\sum_{n\ge 1} \left(
\varphi (\gamma ) e^{-\gamma a} \right)^{a^n},
\end{eqnarray*}
where $\phi(\gamma)={\mathbf E} \exp (\gamma v )$.
First, choose $\gamma >0$ such that $\phi(\gamma)<\infty$.
Then, for any integer $m=1,2,\ldots$, choose $a_m \ge \max ({\mathbf E} v ,2)$ such that
$$
q_m= \varphi (\gamma ) e^{-\gamma a_m} <\frac 1 {2^m}.
$$
So ${\mathbf P} (D(a_m))\le 2^{-m} \to 0$ as $m\to\infty$.

For any $m$ and any $c$,
$$
\Bigl\{\lim\sup_{n\to\infty} \frac{{h}(n)}n > c\Bigr\}
\subseteq
D(a_m) \bigcup \bigl\{\lim\sup_{t\to\infty}
\frac{{h} (n)} n > c\bigr\} \cap \overline{D}(a_m)
$$
and
$$
{\mathbf P}
\left(\bigl\{\lim\sup_{n\to\infty} \frac{{h} (n)} n
> c\bigr\} \cap \overline{D}(a_m)\right) \\
\le \sum_n P(n,c,m)~,
$$
where $
P(n,c,m) =
{\mathbf P}
\left(\bigl\{\frac{{h}(n)} n > c\bigr\} \cap \overline{D}(a_m)\right)$.

We deduce from the union bound that, for all $m$,
\begin{eqnarray*}
P(n,c,m) \le
a_m^n{\mathbf P} \left(\sum_{i=1}^n s_i >cn\right).
\end{eqnarray*}
The inequality follows from the assumption that $v$-family and $s$-family of random
variables are independent.
Hence, by Chernoff's inequality,
\begin{eqnarray*}
P(n,c,m) \le
a_m^n (\psi (\delta))^n  e^{-\delta cn} ,
\end{eqnarray*}
where $\psi (\delta)= {\mathbf E} e^{\delta s}$.
Take $\delta >0$ such that $\psi (\delta )$ is finite and then
$c_m>0$ such that
$$ h_m= a_m\psi (\delta )  e^{-\delta c_m} <1.$$
Then
$$ \sum_{k\in \N}  h_m^k <\infty.$$
Hence for all $m$,
$$ \limsup_n \frac {{h}(n)} n 1_{\overline{D}(a_m)} \le c_m 1_{\overline{D}(a_m)},\quad \mbox{a.s.}$$
Let $\mu$ be a random variable taking the value $c_m$ on the
event
$\overline{D}(a_m)\setminus \overline{D}(a_{m-1})$.
Then $\mu$ is finite a.s. and
$$ \limsup_n \frac {{h}(n)} n \le \mu ,\quad \mbox{a.s.}$$
But $\limsup_n \frac{{h}(n)}{n}$ must be a constant (by ergodicity)
and then this constant is necessarily finite. Indeed,
since
$$\limsup_n \frac{{h}(n)}{n} \ge
 \limsup_n \frac{{h}(n)\circ \theta^{-1}}{n} \quad \mbox{a.s.,}
$$
and since the shift $\theta$ is ergodic, for each $c$,
the event  $\{\limsup_n \frac {{h}(n)} n \le c\}$ has either probability
1 or 0.
\hfill $\Box$

Recall that if $\lambda_c$ is the maximal value of intensity $\lambda$
such that Boolean Model 5 has a.s. finite clumps, for any $\lambda <
\lambda_c$.

\begin{cor}\label{keycor}
Let ${\HH}(t)= \HH_t^0$ be the height at $0\in \Z^d$ in the backward Poisson hail growth model
defined in (\ref{eq:simple2b}).
Under the assumptions of Theorem \ref{LinearGrowth}, for all $\lambda<\lambda_c$, with $\lambda_c>0$
the critical intensity defined above, 
there exists a finite constant $\kappa(\lambda)$ such that
\begin{equation}
\lim\sup_{t\to\infty} \frac{\HH(t)} t = \kappa(\lambda) \quad a.s.
\end{equation}
with $\lambda$ the intensity of the Poisson rain.
\end{cor}

{\sc Proof.}
The proof of the fact that the limit is finite is immediate from 
bound (\ref{eq:bornfond}) and Theorem 2. 
The proof that the limit is constant follows from the ergodicity of the underlying model.
\hfill $\Box$

\begin{lem}\label{leminvar}
Let $a<\lambda_c$, where $\lambda_c$ is the critical value defined above.
For all $\lambda<\lambda_c$,
\begin{equation}\label{eq:invar}
\kappa(\lambda)=\frac{\lambda} a \kappa(a)~,
\end{equation}
\end{lem}

{\sc Proof.}
A Poisson rain of intensity $\lambda$ on the interval $[0,t]$ can be seen
as a Poisson rain of intensity $a$ on the time interval $[0,\lambda t/a]$.
Hence, with obvious notation
$$ {\HH}(t,\lambda)= {\HH}\left(\frac{t\lambda}{a},a\right),$$
which immediately leads to (\ref{eq:invar}).
\hfill $\Box$

\section{Service and Arrivals}\label{secservar}
Below, we focus on the equations for the dynamical system with service and arrivals, namely
on Poisson hail on a hot ground.

Let $W_t^x$ denote the residual workload at $x$ and $t$, namely the
time elapsing between $t$ and the first epoch when the system is free of all workload arrived before time $t$
and intersecting location $x\in \R^d$.
We assume that $H_{0}^x \equiv 0$. Then, with the notation of Section \ref{GM},
\begin{equation}\label{eq:simpledyn}
W_t^x
=
\left(\sigma_{\tau^x(t)}^x
-t+ \tau^x(t) + \sup_{y\in C^x_{\tau^x(t)}} W_{\tau^x(t)}^y
\right)^+
1_{\tau^x(t)\ge 0}.
\end{equation}
We will also consider the Loynes' scheme associated with (\ref{eq:simpledyn}), namely the
random variables
$$ \W_t^x = W_t^x \circ \theta_{t},$$
for all $x\in \R^d$ and $t>0$.
We have
\begin{equation}\label{eq:simpledynb}
\W_t^x
=
\left(\sigma_{\tau^x_-(t)}^x
+ \tau^x_-(t) + \sup_{y\in C^x_{\tau^x_-(t)}} \W_{t+\tau^x_-(t)}^y\circ \theta_{\tau^x_-(t)}
\right)^+
1_{\tau^x_-(t)\ge -t}.
\end{equation}

Assume that $W_0^x=\W_0^x=0$ for all $x$. Using the Loynes-type
arguments (see, e.g., \cite{Loynes} or \cite{Stoyan}), it is easy
to show that for all $x$, $\W_t^x$ is non decreasing in $t$. Let
$$\W_\infty^x=\lim_{t\to \infty} \W_t^x.$$
By a classical ergodic theory argument, the limit $\W_{\infty}^x$ is either
finite a.s. or infinite a.s. Therefore,
 for all integers $n$ and all $(x_1,\ldots,x_n) \in \R^{dn}$,
either $\W_\infty^{x_i}=\infty$ for all $i=1,\ldots,n$ a.s. or
$\W_\infty^{x_i}<\infty$ for all $i=1,\ldots,n$ a.s. In the latter
case,
\begin{itemize}
\item $\{\W_\infty^{x}\}$ is the smallest stationary solution of (\ref{eq:simpledynb});
\item $(\W_t^{x_1},\ldots,\W_t^{x_n})$ converges a.s. to
$(\W_\infty^{x_1},\ldots,\W_\infty^{x_n})$ as $t$ tends to $\infty$.
\end{itemize}

Our main result is (with the notation of Corollary \ref{keycor}):

\begin{theorem}
If $\lambda < \min(\lambda_c,a \kappa(a)^{-1})$,
then for all $x \in \R^{d}$,
$\W_\infty^{x}<\infty$ a.s.
\end{theorem}

{\sc Proof.}
For all $t>0$, we say that $x_0$ is a critical path of length $0$ and span $t$ starting from $x_0$ in the backward growth model
$\{\HH_t^x\}$ defined in (\ref{eq:simple2b}) if $\tau^{x_0}_-(t)<-t$. The height of this path is
$\HH_t^{x_0}=0.$
For all $t>0$, $q\ge 1$, we say that $x_0,x_1,\ldots,x_q$ is a critical path of length $q$ and span $t$ starting from $x_0$
in the backward growth model $\{\HH_t^x\}$ defined in (\ref{eq:simple2b}) if
\begin{eqnarray*}
\HH_t^{x_0} & = & \sigma_{\tau^{x_0}_-(t)}^{x_0} + \HH_{t+\tau^{x_0}_-(t)}^{x_1}\circ \theta_{\tau^{x_0}_-(t)} ~,
\end{eqnarray*}
with $x_1\in C^{x_0}_{\tau^{x_0}_-(t)}$ and $\tau^{x_0}_-(t)>-t$, and if $x_1,\ldots,x_q$ is a critical path of length $q-1$
and span $t+\tau^{x_0}_-(t)$ starting from $x_1$ in the backward growth model
$\{\HH_{t+\tau^{x_0}_-(t)}^x\circ \theta_{\tau^{x_0}_-(t)}\}$. The height of this path
is $\HH_t^{x_0}$.

Assume that $\W_\infty^{x_0}=\infty$. Since $\W_t^x$ is a.s. finite
for all finite $t$ and all $x$, there must exist an increasing sequence $\{t_k\}$, with $t_k\to \infty$,
such that $\W_{t_{k+1}}^{x_0}> \W_{t_{k}}^{x_0}>0$ for all $k$. This in turn implies the
existence, for all $k$, of a critical path of length $q_k$ and span $t_k$, say $x_0,x^k_1,\ldots,x_{q_k}^k$
of height $\HH^{x_0}_{t_k}$ such that
$$
\W^{x_0}_{t_{k+1}} = \HH^{x_0}_{t_k} - t_k >0,
$$
Then
$$
\frac{\HH^{x_0}_{t_k}}{t_k} \ge 1,
$$
for all $k$
and therefore
$$ \kappa (\lambda) \ge \liminf_{k\to \infty} \frac{\HH^{x_0}_{t_k}}{t_k} \ge 1. $$
Using  (\ref{eq:invar}), we get
$$
\kappa(\lambda)= \frac{\lambda}{a} \kappa(a) \ge 1 \quad \mbox{a.s.}
$$
But this contradicts the theorem assumptions.
\hfill $\Box$\\

\begin{rem}
Theorem \ref{thmain} follows from the last theorem and the remarks that precede it.
\end{rem}
\begin{rem}
We will say that the dynamical system with arrivals and service percolates
if there is a time for which the directed graph of RACS present in the system at that
time (where directed edges between two RACS represent the precedence constraints
between them) has an infinite directed component.
The finiteness of the Loynes variable is equivalent to
the non-percolation of this dynamical system.
\end{rem}

\section{Bernoulli Hail on a Hot Grid}\label{secbhhg}

The aim of this section is to discuss discrete versions of the Poisson Hail model.,
namely versions where the server is the grid $\Z^d$ rather than the Euclidean space $\R^d$.
Some specific discrete models were already considered in the analysis of the
Poisson hail model (see e.g. Sections \ref{ssds} and \ref{ssdistim}).
Below, we concentrate on the simplest model, emphasize the main differences
with the continuous case and give a few examples of explicit bounds and evolution equations.

\subsection{Models with Bernoulli Arrivals and Constant Services}\label{1DD}

The state space is $\Z$.
All RACS are pairs of neighbouring points/nodes $\{i,i+1\}$, $i\in \Z$ with service
time 1. In other words, such a RACS requires $1$ unit of time for and simultaneous
service from nodes/servers $i$ and $i+1$. For short, a RACS $\{i,i+1\}$ will be
called ''RACS of type $i$''.

Within each time slot (of size 1), the number of RACS of type $i$ arriving
is a Bernoulli-($p$) random variable. All these variables are mutually independent.
If a RACS of type $i$ and a RACS of type $i+1$ arrive in the same time slot, the FIFO tie
is solved at random (with probability $1/2$). The system is empty at time 0, and
RACS start to arrive from time slot $(0,1)$ on.

\subsubsection{The Growth Model}\label{1000}

{\bf (1)} The Graph ${\cal G}(1)$.\\
We define a precedence graph ${\cal G}(1)$ associated with $p=1$
nodes are all $(i,n)$ pairs where
$i\in \Z$ is a type and $n\in \N =\{1,2,\ldots\}$ is a time.
There are directed edges between certain nodes, some of which are deterministic
and some random. These edges represent precedence constraints: an edge
from $(i,n)$ to $(i',n')$ means that $(i,n)$ ought to be served after $(i',n')$.
Here is the complete list of directed edges:

\begin{enumerate}

\item There is either an edge
$(i,n)\to (i+1,n)$ w.p. 1/2
(exclusive) or an edge $(i+1,n)\to (i,n)$ w.p. 1/2;
we call these random edges {\em spatial};

\item The edges $(i,n) \to (i-1,n-1)$, $(i,n)\to (i,n-1)$, and
$(i,n)\to (i+1,n-1)$ exist for all $i$ and $n\ge 2$; we call these random edges {\em time edges}.

\end{enumerate}

Notice that there are at most five directed edges from each node.
These edges define directed paths: for $x_j = (i_j,n_j)$, $j=1,\ldots ,m$,
the path $x_1 \to x_2 \to \ldots \to x_m$ {\it exists} if (and only if)
all edges along this path exist. All paths in this graph are acyclic.
If a path exists, its {\it length} is
the number of nodes along the path, i.e. $m$.  \\

\noindent
{\bf (2)} The Graph ${\cal G}(p)$.\\
We obtain ${\cal G}(p)$ from ${\cal G}(1)$ by the following thinning:

\begin{enumerate}

\item
Each node of ${\cal G}(1)$ is colored "white" with probability $1-p$ and "black" with
probability $p$, independently of everything else;

\item
If a node is coloured white, then each directed spatial edge from this node
is deleted (recall that there are at most two such edges);

\item

For $n\ge 2$, if a node $(i,n)$ is coloured white, then two time edges
$(i,n) \to (i-1,n-1)$ and $(i,n)\to (i+1,n-1)$ are deleted, and only
the "vertical" one, $(i,n)\to (i,n-1)$, is kept.

\end{enumerate}

So, the sets of nodes are hence the same in ${\cal G}(1)$ and ${\cal G}(p)$ whereas the
set of edges in ${\cal G}(p)$ is a subset of that in ${\cal G}(1)$.
Paths in ${\cal G}(p)$ are defined as above (a path is made of a
sequence of directed edges present in ${\cal G}(p)$).
The graph ${\cal G}(p)$ describes the precedence relationship between RACS
in our basic growth model.

\paragraph{The Monotone Property.}

We have the following monotonicity in $p$:
the smaller $p$, the thinner the graph. In particular,
by using the natural coupling, one can make ${\cal G}(p)\subset {\cal G}(q)$ for
all $p\le q$; here inclusion means that the sets of nodes in both
graphs are the same and the set of edges of ${\cal G}(p)$ is
included in that of ${\cal G}(q)$.\\

\subsubsection{The Heights and The Maximal Height Function}

We now associate {\it heights} to the nodes: the height of a white node is 0
and that of a black one is 1. The {\it height of a path} is the sum of the
heights of the nodes along the path. Clearly, the height of a path
cannot be bigger than its length.

For all $(i,n)$,
let $H_{n}^i=H_{n}^i(p)$ denote the height of the maximal height
path among all paths of ${\cal G}(p)$ which start from node $(i,n)$.
By using the natural coupling alluded to above, we get that
$H_{n}^n(p)$ can be made a.s. increasing in $p$.

Notice that, for all $p\le 1$, for all $n$ and $i$, the random variable $H_{n}^i$ is finite a.s.
To show this, it is enough to consider the case $p=1$ (thanks to monotonicity) and $i=0$
(thanks to translation invariance).
Let
$$
t^+_{n,n} = \min \{ i\ge 1 \ : \ (i,n)\to (i-1,n) \}
$$
and, for $m=n-1,n-2,\ldots ,1$, let
$$
t^+_{m,n} = \min \{ i > t^+_{m+1,n} +1 \ : \ (i,m) \to (i-1,m) \}.
$$
Similarly, let
$$
t^-_{n,n} = \max \{ i\le -1 \ : \ (i,n)\to (i+1,n) \}
$$
and, for $m=n-1,n-2,\ldots ,1$, let
$$
t^-_{m,n} = \max \{ i < t^-_{m+1,n} -1 \ : \ (i,m) \to (i+1,m) \}.
$$
Then all these random variables are finite a.s. (moreover, have finite
exponential moments) and the following rough estimate holds:
$$
H_{n}^0 \le \sum_{i=1}^n \left(
t^+_{1,n} -t^-_{1,n}
\right) +n.
$$

\subsubsection{Time and Space Stationarity}

The driving sequence of RACS is i.i.d. and does not depend on the
random ordering of neighbours which is again i.i.d., so the model is homogeneous
both in time $n=1,2,\ldots$ and
in space $i\in \Z$.
Then we may extend this relation to non-positive indices of $n$ and then introduce
the measure preserving time-transformation
$\theta$ and its iterates $\theta^m, -\infty < m < \infty$. So
$
H_{n}^i \circ \theta^m 
$
is now representing the height of the node $(i,n+m)$ in the model which
starts from the empty state at time $m$.
Again, due to the space homogeneity, for any fixed $n$,
the distribution of the random variable $H_{n}^i$ does not depend on $i$.
So, in what follows, we will write for short
$$
H_n \equiv H_n(p)= H_{n}^i,
$$
when it does not lead to a confusion.

{\bf Definition of function $h$}.
We will also consider paths from $(0,n)$ to $(0,1)$ and we will denote by $h_n=h_n(p)$
the maximal height of all such paths.
Clearly, $h_n \le H_{n}$ a.s.

\subsubsection{Finiteness of the Growth Rate and Its Continuity at 0}

\begin{lem}\label{L1}
There exists a positive probability $p_0\ge 2/5$ such that, for any $p<p_0$,
\begin{equation}\label{Hp}
\limsup_{n\to\infty} H_n/n \le C(p) < \infty, \quad
\mbox{a.s.}
\end{equation}
and
\begin{equation}\label{hp}
h_{n}(p)/n \to \gamma(p) \quad
\mbox{a.s. and in $L_1$},
\end{equation}
with $\gamma (p)$ and $C(p)$
positive and finite constants, $\gamma (p) \le C(p)$.
\end{lem}

\paragraph{Remark.} The sequence $\{ H_n\}$ is neither sub- nor super-additive.

\begin{lem}\label{L11}
For all $p$,
\begin{equation}\label{limH}
\limsup_{n\to\infty} H_n(p)/n \le 2 \gamma (p) \quad \mbox{a.s.}
\end{equation}
\end{lem}

\begin{lem}\label{L2}
Under the foregoing assumptions,
$$
\lim_{p\downarrow 0} \limsup_{n\to\infty} H_n/n = 0
\quad \mbox{a.s.}
$$
\end{lem}

Proofs of Lemmas 5--7 are in a similar spirit to that of the main
results (Borel-Cantelli lemma, branching upper bounds, and also
superadditivity),
and therefore are omitted.

\subsubsection{Exact Evolution Equations for the Growth Model}\label{ee}

We now describe the
exact evolution of the process defined in \S\ref{1000}.
We adopt here the continuous-space interpretation
where a RACS of type $i$ is a segment of length 2 centered in $i\in \Z$.

\begin{figure}[h]
\begin{center}
\includegraphics[width=.8\textwidth]{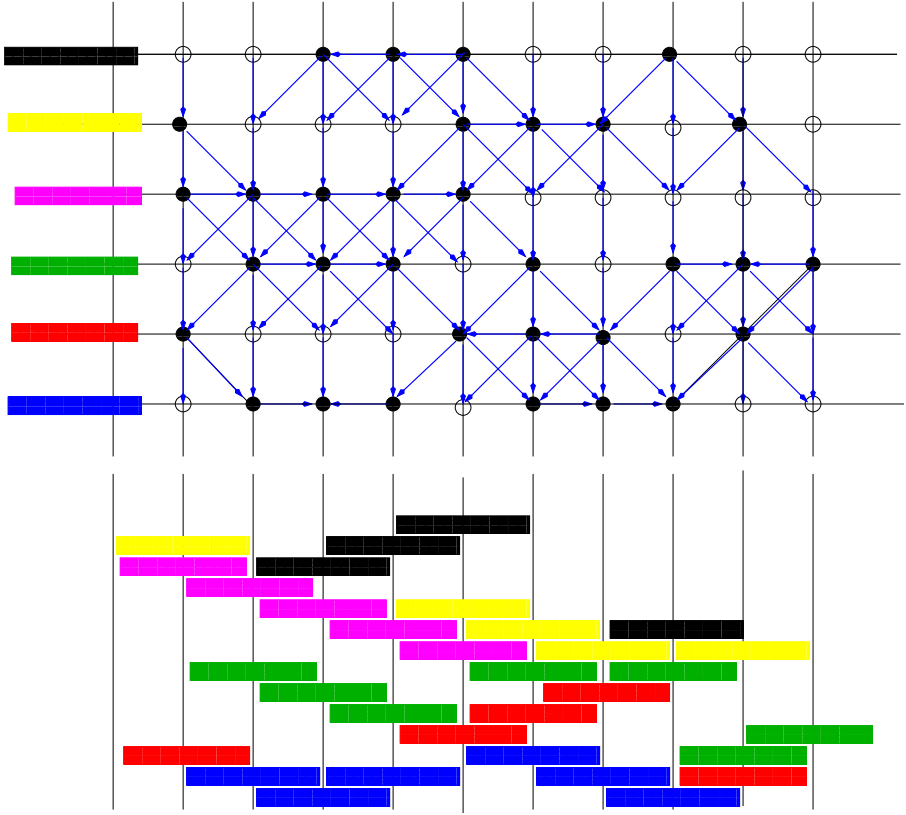}
\end{center}
\caption{Top: A realization of the random graph ${\cal G}(p)$. Only the first 6 time-layers
are represented. A black node at $(i,n)$ represents the arrival of a RACS of type $i$ at time $n$.
Bottom: the associated the heap of RACS, with a visualization of
the height $H_n^i$ of each RACS.}
\label{fig1}
\end{figure}

The variable $H_n^i$ is the height of the {\em last} RACS (segment) of type $i$ that arrived
among the set with time index less than or equal to $n$ (namely with index
$1\le k\le n$), in the growth model under consideration.
If $(i,n)$ is black, then $H_n^i$ is at the same time the height of the maximal
height path starting from node $(i,n)$ in ${\cal G}(p)$
and the height of the RACS $(i,n)$ in the growth model.
If $(i,n)$ is white and the last arrival of type $i$ before time $n$
is $k$, then $H_n^i=H_k^i$. This is depicted in Figure \ref{fig1}

If there are no arrivals of type $i$ in this time interval, then $H_n^i=0$. In general,
if $\beta_n^i$ is the number of segments of type $i$ that arrive in $[1,n]$, then $H_n^i \ge \beta_n^i$.
Let $v_n^i$ be the indicator of the event that $(i,n)$ is an arrival
($v_n^i=1$ if it is black and $v_n^i=0$ otherwise).

Let $e_n^{i,i+1}$ indicate the direction of the edge between $(i,n)$
and $(i+1,n)$: we write $e_n^{i,i+1}=r$ if the right node has priority,
$e_n^{i,i+1}=l$ if the left node has priority.

The following evolution equations hold:
if $v_{n+1}^i=1$, then
\begin{eqnarray*}
H_{n+1}^i &=& (H_{n+1}^{i+1}+1) {\bf I} (e_n^{i,i+1}=r,v_{n+1}^{i+1}=1)\\
&\vee &
(H_{n+1}^{i-1}+1) {\bf I} (e_n^{i-1,i}=l,v_{n+1}^{i-1}=1)\\
&\vee &
(H_n^i
\vee H_n^{i-1}
\vee
H_n^{i+1}
+1)
\end{eqnarray*}
and if $v_{n+1}^i=0$ then $H_{n+1}^i=H_n^i$.
Here, for any event $A$,  ${\bf I}(A)$ is its indicator function:
it equals $1$ if the event occurs and $0$ otherwise.

The evolution equations above may be rewritten as
\begin{eqnarray*}
H_{n+1}^i &=& (H_{n+1}^{i+1}+1) {\bf I} (e_n^{i,i+1}=r,v_{n+1}^{i+1}=1, v_{n+1}^i=1)\\
&\vee &
(H_{n+1}^{i-1}+1) {\bf I} (e_n^{i-1,i}=l,v_{n+1}^{i-1}=1, v_{n+1}^i=1)\\
&\vee &
(H_n^i \vee H_n^{i-1}{\bf I} (
v_{n+1}^i=1)
\vee H_n^{i+1}{\bf I} (
v_{n+1}^i=1)+
{\bf I} (v_{n+1}^i=1)).
\end{eqnarray*}

\subsubsection{Exact Evolution Equations for the Model with Service}\label{eed}

The system with service can be described as follows: there is an infinite number of servers, each of which
serves with a unit rate. The servers are located at points
$1/2+i$, $-\infty < i < \infty$. For each $i$,  RACS $(i,n)$
(or customer $(i,n)$)
is a customer of ``type'' $i$ that arrives with probability $p$
at time $n$ and needs one unit of time
for simultaneous service from two servers located at points
$i-1/2$ and $i+1/2$. So, at most one customer of each type arrives
at each integer time instant. If customers of types $i$ and $i+1$
arrive at time $n$, then one makes a decision, that either
$i$ arrives earlier or $i+1$ arrives earlier, 
at random with equal probabilities,
$$
{\mathbf P} (\mbox{customer}
\ \ i \ \ \mbox{arrives earlier than customer} \ \ i+1) =
{\mathbf P} (e_{n}^{i,i+1}= l) = 1/2.
$$
Each server serves customers
in the order of arrival. A customer leaves the system after the
completion
of its service.
As before, we may assume that, for each $(i,n)$, customer $(i,n)$
arrives with probability $1$, but is ``real''(``black'') with
probability
$p$ and ``virtual''(``white'') with probability $1-p$.

Assume that the system is empty at time $0$ and that the first customers arrive
at time $1$. Then, for any $n=1,2,\ldots$, the quantity
$W_{n}^{i}:=\max (T_{n}^{i}-(n-1),0)$ is the residual amount of time (starting
from time $n$) which is needed for the last
real customer of type $i$ (among customers $(i,1),\ldots , (i,n)$)
to receive the service (or equals zero if there is no real customers there).

Then these random variables satisfy the equations, for $n\ge 1, -\infty < i < \infty$,
\begin{eqnarray*}
W_{n+1}^i &=& (W_{n+1}^{i+1}+1) {\bf I} (e_n^{i,i+1}=r,v_{n+1}^{i+1}=1, v_{n+1}^i=1)\\
&\vee &
(W_{n+1}^{i-1}+1) {\bf I} (e_n^{i-1,i}=l,v_{n+1}^{i-1}=1, v_{n+1}^i=1)\\
&\vee &
((W_n^i-1)^+ + {\bf I} (v_{n+1}^i=1))\\
&\vee &
((W_n^{i-1}-1)^+ +1) {\bf I} (v_{n+1}^i=1)\\
&\vee &
((W_n^{i+1}-1)^+ +1) {\bf I} (v_{n+1}^i=1).
\end{eqnarray*}
Since the heights are equal to 1 (and time intervals have length 1),
the last two terms in the equation may be simplified, for instance,
$
((W_n^{i-1}-1)^+ +1) {\bf I} (v_{n+1}^i=1)
$
may be replaced by
$W_n^{i-1} {\bf I} (v_{n+1}^i=1).$

In the case of random heights $\{\sigma_n^i\}$,
the random variables $\{ W_n^i\}$ satisfy the recursions
\begin{eqnarray*}
W_{n+1}^i &=& (W_{n+1}^{i+1}+\sigma_{n+1}^i) {\bf I} (e_n^{i,i+1}=r,v_{n+1}^{i+1}=1, v_{n+1}^i=1)\\
&\vee &
(W_{n+1}^{i-1}+\sigma_{n+1}^i) {\bf I} (e_n^{i-1,i}=l,v_{n+1}^{i-1}=1, v_{n+1}^i=1)\\
&\vee &
((W_n^i-1)^+ + \sigma_{n+1}^i{\bf I} (v_{n+1}^i=1))\\
&\vee &
((W_n^{i-1}-1)^+ +\sigma_{n+1}^i) {\bf I} (v_{n+1}^i=1)\\
&\vee &
((W_n^{i+1}-1)^+ +\sigma_{n+1}^i) {\bf I} (v_{n+1}^i=1).
\end{eqnarray*}

The following monotonicity property holds: for any $n$ and $i$,
$$
W_{n+1}^{i} \circ \theta^{-n-1} \le W_{n}^{i} \circ \theta^{-n} \quad
\mbox{a.s.}
$$
Let
$$
p_{0} = \sup \{ p \ : \ \Gamma (p) \le 1 \}.
$$

\begin{theorem}
If $p<p_{0}$, then, for any $i$, random variables
$W_{n}^{i}$ converge weakly to a
proper limit. Moreover, there exists a stationary random
vector $\{ W^{i}, -\infty < i < \infty \}$ such that, for any
finite integers $i_{0}\le 0 \le i_{1}$, the finite-dimensional random
vectors
$$
(W_{n}^{i_{0}}, W_{n}^{i_{0}+1}, \ldots , W_{n}^{i_{1}-1},
W_n^{i_{1}})
$$
converge weakly to the vector
$$
(W^{i_{0}}, W^{i_{0}+1}, \ldots , W^{i_{1}-1}.
W^{i_{1}})
$$
\end{theorem}

\begin{theorem}
If $p<p_{0}$, then the random variables
$$
\min \{ i\ge 0 \ : \  W^{i}=0 \}
\quad
\mbox{and}
\quad
\max \{ i \le 0 \ : \  W^{i}=0 \}
$$
are finite a.s.
\end{theorem}

\section{Conclusion}
We conclude with a few open questions. The first class of questions pertain to stochastic
geometry \cite{Stoyan:1995:SGA}:
\begin{itemize}
\item How does the RACS exclusion process which is that of the RACS in service at time $t$ in
steady state compare to other exclusion processes (e.g. Mat\'ern, Gibbs)?
\item Assuming that the system is stable,
can the {\em undirected} graph of RACS present in the steady state regime
percolate?
\end{itemize}
The second class of questions are classical in queueing theory and pertain to existence and properties
of the stationary regime:
\begin{itemize}
\item In the stable case, does the stationary solution $\W_{\infty}^0$
always have a light tail? At the moment, we can show this under
extra assumptions only. Notice that in spite of the fact that the Poisson hail model falls in the
category of infinite dimensional max plus linear systems.
Unfortunately, the techniques developed for analyzing the tails of
the stationary regimes of finite dimensional max plus linear
systems \cite{moments} cannot be applied here.
\item In the stable case, does the Poisson hail equation (\ref{eq:simpledynb}) admit other stationary regimes than obtained from
$\{\W_{\infty}^x\}_x$, the minimal stationary regime?
\item For what other service disciplines still respecting the hard exclusion rule
like e.g. priorities or first/best fit can one also construct a steady state?
\end{itemize}

\section{Appendix}\label{GetInd}

\begin{prop}\label{Sergey}
For any pair $(X,Y)$ of random variables with light-tailed
marginal distributions, there exists a coupling with another pair
$(\xi,\eta)$ of i.i.d. random variables with a common light-tailed
marginal distribution and such that
$$
\max (X,Y) \le \min (\xi,\eta) \quad \mbox{a.s}.
$$
\end{prop}

{\sc Proof.} Let $F_X$ be the distribution function of $X$ and $F_Y$
the distribution function of $Y$. Let $C>0$ be such that ${\mathbf
E} e^{CX}$ and ${\mathbf E} e^{CY}$ are finite. Let $\zeta = \max
(0,X,Y)$. Since $e^{C\zeta}\le 1+e^{CX}+e^{CY}$, $\zeta$ also has
a light-tailed distribution, say $F$.

Let $\overline{F}(x)=1-F(x)$,
$\overline{G}(x)=\overline{F}^{1/2}(x)$, and
$G(x)=1-\overline{G}(x)$. Let $\xi$ and $\eta$ be i.i.d. with
common distribution $G$. Then ${\mathbf E}^{c\xi}$ is finite, for
any $c<C/2$.

Finally, a coupling of $X,Y,\xi$, and $\eta$ may be built as
follows. Let $U_1$, $U_2$ be two i.i.d. random variable having
uniform $(0,1)$ distribution. Then let $\xi = G^{-1}(U_1)$, $\eta =
G^{-1}(U_2)$ and $\zeta = \min (\xi , \eta )$. Finally, define $X$
and $Y$ given $\max (X,Y)= \zeta$ and conditionally independent of
$(\xi,\eta)$.

\hfill $\Box$

\end{document}